\documentclass[a4paper,11pt]{amsart}
\usepackage{amssymb,amsfonts,amsxtra,amscd,mathrsfs, stmaryrd,bbm}
\usepackage[all]{xy}
\usepackage{fullpage}
\theoremstyle{plain}
\newtheorem{theorem}{Theorem}[section]
\newtheorem{lemma}[theorem]{Lemma}
\newtheorem{cor}[theorem]{Corollary}
\newtheorem{prop}[theorem]{Proposition}
\theoremstyle{definition}
\newtheorem{defi}[theorem]{Definition}

\theoremstyle{remark}
\newtheorem{rem}[theorem]{Remark}
\numberwithin{equation}{section}
\newcommand{\ai}{\ensuremath{A_\infty}}
\newcommand{\forms}[2][\bullet]{\ensuremath{\Omega^{#1}(#2)}}
\newcommand{\cforms}[2][\bullet]{\ensuremath{\widehat{\Omega}{\vphantom{\Omega}}^{#1}(#2)}}
\newcommand{\hlie}{\ensuremath{\mathfrak{h}_{2n|m}}}
\newcommand{\glie}{\ensuremath{\mathfrak{g}_{2n|m}}}
\newcommand{\tglie}{\ensuremath{\tilde{\mathfrak{g}}_{2n|m}}}
\newcommand{\hl}[1]{\ensuremath{\mathfrak{h}[#1]}}
\newcommand{\gl}[1]{\ensuremath{\mathfrak{g}[#1]}}
\newcommand{\tgl}[1]{\ensuremath{\tilde{\mathfrak{g}}[#1]}}
\newcommand{\osp}{\ensuremath{\mathfrak{osp}_{2n|m}}}
\newcommand{\gc}[1][\bullet]{\ensuremath{\mathcal{G}_{#1}}}
\newcommand{\cek}[1][\bullet]{\ensuremath{C_{#1}(\glie,\osp)}}
\newcommand{\stlim}[1][\bullet]{\ensuremath{C_{#1}(\mathfrak{g},\mathfrak{osp})}}
\newcommand{\chd}[1]{\ensuremath{\mathscr{C}(#1)}}
\newcommand{\grpchd}[1]{\ensuremath{\Gamma_{k_1,\ldots,k_l}(#1)}}
\newcommand{\csalg}[1]{\ensuremath{\widehat{S}(#1)}}
\newcommand{\salg}[1]{\ensuremath{S(#1)}}
\newcommand{\gf}{\ensuremath{\mathbb{R}}}
\newcommand{\id}{\ensuremath{\mathbbm{1}}}
\newcommand{\dilim}[2]{\ensuremath{\varinjlim_{#1} #2}}
\newcommand{\innprod}{\ensuremath{\langle -,- \rangle}}
\newcommand{\pder}[1]{\ensuremath{\partial_{#1}}}
\newcommand{\noproof}{\begin{flushright} \ensuremath{\square} \end{flushright}}
\DeclareMathOperator{\Hom}{Hom}
\DeclareMathOperator{\Der}{Der}
\DeclareMathOperator{\Coder}{Coder}
\DeclareMathOperator{\tr}{tr}

\begin{document}
\begin{abstract}
We give a conceptual formulation of Kontsevich's `dual construction' producing graph cohomology classes from a differential graded Frobenius algebra with an odd scalar product. Our construction -- whilst equivalent to the original one -- is combinatorics-free and is based on the Batalin-Vilkovisky formalism, from which its gauge-independence is immediate.
\end{abstract}
\title{Graph cohomology classes in the Batalin-Vilkovisky formalism}
\author{Alastair Hamilton}
\address{Max-Planck-Institut f\"ur Mathematik, Vivatsgasse 7, 53111 Bonn, Germany.}
\email{hamilton@mpim-bonn.mpg.de}
\author{Andrey Lazarev}
\address{IH\'ES, Le Bois-Marie, 35 route de Chartres, F-91440, Bures-sur-Yvette,
France.}
\email{lazarev@ihes.fr}
\maketitle
\tableofcontents
%
%
%

\section{Introduction}

M. Kontsevich introduced graph complexes in the beginning of the nineties \cite{kontsg}, \cite{kontfd} and proposed two constructions producing graph homology and cohomology classes. The first construction uses an $\ai$-algebra with an invariant scalar product (also known as a \emph{symplectic} or a \emph{cyclic} $\ai$-algebra \cite{htopnc}) as the input and gives rise to a family of homology classes in the ribbon graph complex. One can also generalise this construction to algebras over other cyclic operads (such as cyclic $L_\infty$ or $C_\infty$-algebras) to produce classes in the appropriate versions of the graph complex considered in, e.g. \cite{CV} and \cite{GK}. This so called `direct construction' is by now well-understood from the point of view of the homology of Lie algebras. Topological conformal field theories (TCFTs) provide an alternative conceptual explanation of the direct construction and we refer the reader to \cite{ccs} for a detailed discussion of these approaches.

The `dual construction' starts from a contractible differential graded Frobenius algebra with an odd scalar product and produces a family of \emph{cohomology} classes in the ribbon graph complex. It draws its motivation from the quantum Chern-Simons theory, cf. \cite{AxS1}, \cite{AxS2}, \cite{bar} and \cite{kontfd}; the aforementioned contractible differential graded Frobenius algebra being analogous to the de Rham algebra on a three-dimensional sphere.

In contrast with the direct construction, the dual one is not at all well-understood. It is unclear, first of all, why it gives rise to a graph cocycle. It is then unclear why the class of this cocycle is independent of the choices involved (gauge-independence). It turns out, in fact, that one has to impose a rather stringent condition on the Frobenius algebra for it to give a cocycle; otherwise the cocycle is obtained on a certain partial compactification of the moduli space of ribbon graphs (or the corresponding compactification of the moduli space of Riemann surfaces). This will be discussed in the joint work of the second author with J. Chuang, as well as an interpretation of the dual construction as a kind of TCFT.

As mentioned above there are many versions of the graph complex of which the most important ones are the ribbon graph complex, the Lie graph complex and the commutative (undecorated) graph complex. The latter is the simplest of all; it is of interest mainly because of its relation to the invariants of homology spheres (and more generally smooth manifolds), cf. for example \cite{BC}.

In the present paper we give a conceptual and definitive formulation of the dual construction from the point of view of homological algebra, analogous to our treatment of the direct construction in \cite{ccs}. We only treat the case of the \emph{commutative} graph complex; more work needs to be done to accommodate other flavours of graph cohomology. This is essentially because the graphs under consideration do not have any \emph{loops}, the standard source of difficulties in quantum field theory.

Central to our approach is the Batalin-Vilkovisky (BV) formalism, devised by Russian physicists Batalin and Vilkovisky \cite{BV} as a method of quantising gauge theories. More recently the BV-formalism was applied to the problem of deformation quantisation, cf. \cite{CF}. A modern formulation of the BV-formalism was given in the fundamental paper by A. Schwarz \cite{schwarz}.

Our formulation utilises a theorem of Kontsevich, which expresses graph homology through the homology of a certain infinite-dimensional Lie algebra $\mathfrak{g}$ consisting of Hamiltonian vector fields vanishing at the origin. Very roughly, for a given Hamiltonian $h\in \mathfrak{g}$ and a contractible differential graded commutative Frobenius algebra $A$, we construct a certain Hamiltonian $\Psi_A$ on an odd symplectic superspace; integrating this superfunction over an appropriate Lagrangian subspace gives a Chevalley-Eilenberg $1$-cochain on $\mathfrak{g}$. This construction then generalises to give cochains in all dimensions.

It turns out that the constructed cochain is a cocycle; this is highly nontrivial and follows from the Batalin-Vilkovisky yoga together with the fact that the vector field generated by the Hamiltonian $\Psi_A$ is divergence-free; this fact is in turn related to the absence of loops in the commutative graph complex. Using physical language one can say that our model has no quantum anomalies. In this connection we mention the paper of A. Schwarz \cite{Sch2} where similar ideas were spelled out. Furthermore, the BV-formalism also implies the independence of the obtained cohomology class on the choice of a gauge-condition, i.e. the choice of the Lagrangian subspace.

The paper is organised as follows. In Section \ref{sec_extcalc} we recall the basics of symplectic geometry on a flat superspace and introduce the Lie algebras of Hamiltonian vector fields. This material is fairly standard and is treated in detail in, e.g. \cite{htopnc}. Section \ref{sec_grphom} is a recollection of graph homology and Kontsevich's theorem relating it to the homology of a certain infinite-dimensional Lie algebra. Here our exposition follows \cite{hamgraph}. Section \ref{sec_wick} treats the Wick rotation in the setting of finite-dimensional superspaces. This material seems to be completely standard for physicists but we have decided to include it in order to make the paper self-contained. Section \ref{sec_BVform} gives the basics of the BV-geometry. There are many good papers where the BV-geometry is considered from various standpoints, besides Schwarz's paper mentioned earlier; particularly \cite{Get}, \cite{Kh} and \cite{KS}. We don't explicitly use the more general results of these papers. Instead we give a simple and direct treatment of the BV-formalism for flat spaces which is sufficient for our purposes. Finally, in Section \ref{sec_dual} we give the BV-formulation of the dual construction and show its equivalence to the original combinatorial version.

\subsection{Notation and conventions}

Throughout the paper we shall work with vector superspaces over the real numbers. A vector superspace $W:=W_0\oplus W_1$ is a $\mathbb{Z}/2\mathbb{Z}$-graded vector space. Our tendency will be to omit the adjective `super' wherever possible. Given a vector (super)space $W$, its parity reversion $\Pi W$ is defined by the formulae:
\[ \Pi W_0:=W_1, \quad \Pi W_1:=W_0. \]
The dual space of $W$ will be denoted by $W^*$. The (super)trace of a linear operator $f:W\to W$ will be denoted by $\tr(f)$.

From the vector space $W$, we may form the tensor algebra $T(W)$ and we denote the Lie algebra of derivations by $\Der[T(W)]$. Likewise, we may form the tensor \emph{coalgebra} $T(W)$ (with the noncocommutative comultiplication) and we denote the Lie algebra of \emph{coderivations} by $\Coder[T(W)]$. Our convention will be to refer to elements of $\Der[T(W)]$ as \emph{vector fields}. We have the following correspondence between multilinear maps and coderivations:
\begin{equation} \label{eqn_correspondence}
\begin{array}{ccc}
\Coder[T(W)] & \cong & \Hom(T(W),W), \\
\xi & \mapsto & p\circ\xi;
\end{array}
\end{equation}
where $p:T(W)\to W$ is the canonical projection. We denote the Lie subalgebra of $\Coder[T(W)]$ consisting of coderivations which vanish on $\gf$ by $\Coder_+[T(W)]$. Likewise, we denote the Lie algebra consisting of vector fields which vanish at zero by $\Der_+[T(W)]$.

The symmetric group on $n$ letters will be denoted by $S_n$. We may similarly define the Lie algebra of (co)derivations for the symmetric (co)algebra $\salg{W}:=\bigoplus_{n=0}^\infty [W^{\otimes n}]_{S_n}$. The correspondence defined by \eqref{eqn_correspondence} continues to hold when we replace $T(W)$ with $\salg{W}$.

For any vector space $W$, we can define mutually inverse maps
\[ i:\salg{W}\to\bigoplus_{n=0}^\infty [W^{\otimes n}]^{S_n}, \quad \pi:\bigoplus_{n=0}^\infty [W^{\otimes n}]^{S_n}\to\salg{W}\]
between $S_n$-coinvariants and $S_n$-invariants by the formulae
\begin{equation} \label{eqn_ipi}
i_n(x_1\otimes\ldots\otimes x_n):=\sum_{\sigma\in S_n} \sigma\cdot[x_1\otimes\ldots\otimes x_n], \quad \pi_n(x_1\otimes\ldots\otimes x_n):=\frac{1}{n!} x_1\otimes\ldots\otimes x_n.
\end{equation}

The initial data for our constructions will be a \emph{differential graded Frobenius algebra}. This is a differential $\mathbb{Z}/2\mathbb{Z}$-graded associative algebra $A$ with a nondegenerate scalar product $\innprod$ satisfying the following conditions:
\begin{align}
\label{eqn_invinnprod} & \langle ab,c \rangle = \langle a,bc \rangle, \quad \text{for all } a,b,c\in A;\\
\label{eqn_dinv} & \langle d(a),b \rangle + (-1)^a\langle a,d(b) \rangle = 0, \quad \text{for all } a,b\in A.
\end{align}

Let $k$ be a positive integer. A partition $c$ of $\{1,\ldots, 2k\}$ such that every $x \in c$ is a set consisting of precisely two elements will be called a \emph{chord diagram}. The set of all such chord diagrams will be denoted by $\chd{k}$.

\section{Exterior calculus on a superspace} \label{sec_extcalc}

In this section we recall the basic geometrical apparatus and terminology used to describe the geometry of (flat) superspaces. The Lie algebras $\glie$ of Hamiltonian vector fields are defined and the definition of Chevalley-Eilenberg homology is recalled.

\subsection{Polynomial forms} \label{sec_commgeom}

We begin with the definition of polynomial valued differential forms.

\begin{defi}
Let $W$ be a vector space, the module of polynomial 1-forms $\forms[1]{W}$ is defined as the quotient of $\salg{W^*}\otimes\salg{W^*}$ by the relations
\[ x \otimes yz=(-1)^{|z|(|y|+|x|)}zx\otimes y + xy\otimes z, \]
for $x,y,z\in \salg{W^*}$.

This is a left $\salg{W^*}$-module via the action
\[ a \cdot (x\otimes y):=ax\otimes y, \]
for $a,x,y\in\salg{W^*}$.

Let $d:\salg{W^*}\to\forms[1]{W}$ be the map given by the formula,
\[ d(x):= 1\otimes x. \]
The map $d$ is a derivation of degree zero.
\end{defi}

\begin{prop}
Let $W$ be a vector space, then the following vector spaces are isomorphic:
\begin{displaymath}
\begin{array}{ccc}
W^*\otimes \salg{W^*} & \cong & \forms[1]{W} \\
x\otimes y & \mapsto & dx\cdot y \\
\end{array}
\end{displaymath}
\end{prop}
\noproof

\begin{defi}
Let $W$ be a vector space and let $A:=\salg{W^*}$. The module of polynomial forms $\forms{W}$ is defined as
\[ \forms{W}:=S_A\left[\Pi\Omega^1(W)\right]=A \oplus \bigoplus_{i=1}^\infty (\underbrace{\Pi\Omega^1(W)\underset{A}{\otimes} \ldots \underset{A}{\otimes} \Pi\Omega^1(W)}_{i \text{ factors}})_{S_i}. \]

Since $\Omega^1(W)$ is a module over the commutative algebra $A$, $\forms{W}$ has the structure of a commutative algebra whose multiplication is the standard multiplication on the symmetric algebra $S_A\left[\Pi\Omega^1(W)\right]$. The map $d:\salg{W^*} \to \forms[1]{W}$ lifts uniquely to a map $d:\forms{W} \to \forms{W}$ which gives $\forms{W}$ the structure of a differential graded commutative algebra.
\end{defi}

\begin{defi} \label{def_diffoperators}
Let $W$ be a vector space and let $\eta:\salg{W^*} \to \salg{W^*}$ be a vector field:
\begin{enumerate}
\item
We can define a vector field $L_\eta: \forms{W} \to \forms{W}$, called the Lie derivative, by the formulae:
\begin{align*}
L_\eta(x)&:=\eta(x), \\
L_\eta(dx)&:=(-1)^{\eta}d(\eta(x)); \\
\end{align*}
for any $x \in \salg{W^*}$.
\item
We can define a vector field $i_\eta:\forms{W} \to \forms{W}$, called the contraction operator, by the formulae:
\begin{align*}
i_\eta(x)&:=0, \\
i_\eta(dx)&:=\eta(x); \\
\end{align*}
for any $x \in \salg{W^*}$.
\end{enumerate}
\end{defi}

These maps can be shown to satisfy the following identities:

\begin{lemma} \label{lem_operatorids}
Let $W$ be a vector space and let $\eta,\gamma:\salg{W^*} \to \salg{W^*}$ be vector fields, then
\begin{enumerate}
\item \label{item_operatoridsone} $L_\eta=[i_\eta,d]$.
\item \label{item_operatoridstwo} $[L_\eta,i_\gamma]=i_{[\eta,\gamma]}.$
\item $L_{[\eta,\gamma]}=[L_\eta,L_\gamma].$
\item \label{item_operatoridsfour} $[i_\eta,i_\gamma]=0.$
\item $[L_\eta,d]=0.$
\end{enumerate}
\end{lemma}

\subsection{Formal geometry}

The geometry described in the last subsection admits a formal analogue, where polynomials are replaced with formal power series. Since the vector spaces we are working with are assumed to be finite-dimensional, this formal geometry is defined by simply taking an adic completion of the definitions given in subsection \ref{sec_commgeom}. We summarise the differences as follows:

For a finite-dimensional vector space $W$ the symmetric algebra $\salg{W^*}$ gets replaced with the \emph{completed symmetric algebra}
\[ \csalg{W^*}:=\prod_{i=0}^\infty S^i(W^*).\]
$\salg{W^*}$ is properly contained inside $\csalg{W^*}$ as a subalgebra.
A \emph{vector field} on $\csalg{W^*}$ is a derivation
\[ \eta:\csalg{W^*}\to\csalg{W^*} \]
which is \emph{continuous} with respect to the adic topology. We denote the Lie algebra of all such vector fields by $\Der[\csalg{W^*}]$. Since $W$ is finite-dimensional, any vector field $\eta:\salg{W^*}\to\salg{W^*}$ extends uniquely to a vector field $\widehat{\eta}:\csalg{W^*}\to\csalg{W^*}$; hence $\Der[\salg{W^*}]$ is properly contained inside $\Der[\csalg{W^*}]$ as a Lie subalgebra.

There exists a formal version of the algebra of forms on $W$ which is denoted by $\cforms{W}$. This differential graded algebra is just the adic completion of $\forms{W}$. The formal versions of the Lie derivative and contraction operators introduced in Definition \ref{def_diffoperators} are defined analogously.

The identity
\[ [\salg{W}]^* = \prod_{i=0}^\infty [(W^*)^{\otimes i}]^{S_i}\]
and the maps $i$ and $\pi$ defined by equation \eqref{eqn_ipi} give rise to the following isomorphism of Lie algebras:
\begin{equation} \label{eqn_coderder}
\begin{array}{ccc}
\Coder[\salg{W}] & \cong & \Der[\csalg{W^*}], \\
\eta & \mapsto & \eta^{\vee}:=\pi\circ\eta^*\circ i;
\end{array}
\end{equation}
which also holds when we replace $\Coder$ and $\Der$ by $\Coder_+$ and $\Der_+$ respectively.

\subsection{Flat symplectic geometry} \label{sec_sympgeom}

In this subsection we recall the basic definitions for symplectic geometry on a flat superspace according to \cite{kontsg}. We begin by giving the definition of a symplectic form. In this paper we will mostly be interested in \emph{constant} 2-forms; these are 2-forms $\omega\in\forms[2]{W}$ which can be written in the form $\omega=\sum_i dy_i dz_i$ for some functions $y_i,z_i\in W^*$. Note that any constant 2-form is closed; that is to say that $d\omega=0$.

\begin{defi} \label{def_symplectic_form}
Let $W$ be a vector space and $\omega \in \forms[2]{W}$ be any \emph{constant} 2-form. We say that $\omega$ is a \emph{symplectic form} if it is nondegenerate, that is to say that the following map is bijective;
\begin{equation} \label{eqn_nondegenerate}
\begin{array}{ccc}
\Der[\salg{W^*}] & \to & \forms[1]{W}, \\
\eta & \mapsto & i_\eta(\omega).
\end{array}
\end{equation}
\end{defi}

\begin{defi}
Let $W$ be a vector space and let $\omega \in \forms[2]{W}$ be a symplectic form. We say a vector field $\eta:\salg{W^*} \to \salg{W^*}$ is a \emph{symplectic vector field} if $L_\eta(\omega)=0$.
\end{defi}

Recall that a (skew)-symmetric bilinear form $\innprod$ on $W$ is called \emph{nondegenerate} if the map
\begin{displaymath}
\begin{array}{ccc}
W & \to & W^* \\
x & \mapsto & [a \mapsto \langle x,a \rangle]
\end{array}
\end{displaymath}
is bijective. The following proposition provides an interpretation of constant 2-forms.

\begin{prop}
Let $W$ be a vector space. There exists a map
\[ \Upsilon:\{ \omega \in \forms[2]{W} : \omega \text{ is constant} \} \to (\Lambda^2 W)^* \]
defined by the formula
\[ \Upsilon(dxdy):=(-1)^x[x\otimes y-(-1)^{xy}y\otimes x] \]
which provides an isomorphism between constant 2-forms and \emph{skew-symmetric} bilinear forms on $W$. Furthermore, a constant 2-form $\omega\in\forms[2]{W}$ is nondegenerate if and only if the corresponding bilinear form $\innprod:=\Upsilon(\omega)$ is nondegenerate.
\end{prop}
\noproof

The analogous definitions for \emph{symmetric} bilinear forms are as follows.

\begin{defi} \label{def_symmforms}
Let $W$ be a vector space. To any quadratic function $\sigma\in S^2(W^*)$ there corresponds a \emph{symmetric} bilinear form
\[ \innprod:W\otimes W\to\gf \]
defined by the formula $\innprod:= i_2(\sigma)$. We say that $\sigma$ is nondegenerate if and only if the corresponding bilinear form $\innprod$ is nondegenerate.
\end{defi}

Since any nondegenerate bilinear form on a vector space provides an identification of that vector space with its dual, it gives rise to a nondegenerate bilinear form on the dual space which is described as follows:

\begin{defi}
Let $\innprod$ be a nondegenerate (skew)-symmetric bilinear form on a vector space $W$, then it gives rise to isomorphisms
\[ D_l:W\to W^* \quad\text{and}\quad D_r:W\to W^* \]
defined by the formulae
\[ D_l(u):=[x\mapsto\langle u,x \rangle] \quad\text{and}\quad D_r(u):=[x\mapsto\langle x,u \rangle].\]
The \emph{inverse form}, denoted by $\innprod^{-1}$, is given by the following commutative diagram:
\begin{equation} \label{fig_inverseform}
\xymatrix{ & \gf \\ W\otimes W \ar[ur]^{\innprod} \ar[rr]^{D_l \otimes D_r} && W^*\otimes W^* \ar[ul]_{\innprod^{-1}}}
\end{equation}
\end{defi}

A symplectic vector space is defined as a vector space $W$ which comes equipped with a symplectic form $\omega\in\forms[2]{W}$. There is a canonical symplectic vector space $V$ of dimension $2n|m$ with an \emph{even} symplectic form. There are canonical coordinates
\begin{equation} \label{eqn_cansymp}
\underbrace{p_1,\ldots,p_n;q_1,\ldots,q_n}_{\text{even}};\underbrace{x_1,\ldots,x_m}_{\text{odd}}
\end{equation}
on $V$ which form a basis for $V^*$ and the symplectic form $\omega\in\forms[2]{V}$ is given by
\[\omega:=\sum_{i=1}^n dp_i dq_i +\frac{1}{2}\sum_{i=1}^m dx_i dx_i.\]

\begin{lemma} \label{lem_fieldsforms}
Let $(W,\omega)$ be a symplectic vector space, then the map
\[ \Phi:\Der[\salg{W^*}] \to \forms[1]{W} \]
given by the formula $\Phi(\eta):=i_\eta(\omega)$ induces a one-to-one correspondence between \emph{symplectic} vector fields and \emph{closed} 1-forms:
\[ \Phi: \{ \eta \in \Der[\salg{W^*}] : L_\eta(\omega) = 0 \} \overset{\cong}{\longrightarrow} \{ \lambda \in \forms[1]{W} : d\lambda = 0 \}. \]
\end{lemma}
\noproof

This means that to any function in $\forms[0]{W}:=\salg{W^*}$ we can associate a certain symplectic vector field;
\begin{displaymath}
\begin{array}{ccc}
\forms[0]{W} & \to & \{ \eta \in \Der[\salg{W^*}] : L_\eta(\omega) = 0 \}, \\
a & \mapsto & \Phi^{-1}(da).
\end{array}
\end{displaymath}
We will typically denote the resulting vector field by the corresponding Greek letter, for instance in this case we have $\alpha:=\Phi^{-1}(da)$.

This allows us to define a Lie algebra as follows.

\begin{defi}
Given a symplectic vector space $(W,\omega)$ with an even symplectic form we can define a Lie algebra $\hl{W}$ whose underlying vector space is $\hl{W}:=S(W^*)$.

The Lie bracket on $\hl{W}$ is given by the formula;
\begin{equation} \label{eqn_bracket}
\{ a,b \}:= (-1)^a L_\alpha(b), \quad \text{for any } a,b \in \hl{W}.
\end{equation}

We define a new Lie algebra $\gl{W}$ as the Lie subalgebra of $\hl{W}$ which has the underlying vector space
\[\gl{W}:=\bigoplus_{i=2}^\infty S^i(W^*).\]
The Lie algebra $\tgl{W}$ is defined to be the Lie subalgebra of $\gl{W}$ which has the underlying vector space
\[\tgl{W}:=\bigoplus_{i=3}^\infty S^i(W^*).\]
\end{defi}

\begin{rem}
The anti-symmetry and Jacobi identity follow straightforwardly from the identities in Lemma \ref{lem_operatorids}, see \cite{ccs} and also Lemma \ref{lem_lieids}. Note that the case when $W$ is a symplectic vector space with an \emph{odd} symplectic form admits an analogous treatment with appropriate modifications to the signs and grading, however its description will be deferred to section \ref{sec_BVform} on the BV-formalism.
\end{rem}

The Lie algebras corresponding to the canonical $2n|m$-dimensional symplectic vector space $V$, cf. \eqref{eqn_cansymp}, are denoted as follows:
\[ \hlie:=\hl{V}, \quad \glie:=\gl{V}, \quad \tglie:=\tgl{V}. \]
The quadratic part of $\glie$ forms a Lie subalgebra of $\glie$ which may be identified with the Lie algebra $\osp$ of linear symplectic vector fields on $V$, hence $\glie$ decomposes as
\[ \glie = \osp\ltimes \tglie. \]
We denote the direct limits of these Lie algebras by
\[ \mathfrak{h}:=\dilim{n,m}{\hlie}, \quad \mathfrak{g}:=\dilim{n,m}{\glie}, \quad \tilde{\mathfrak{g}}:=\dilim{n,m}{\tglie} \quad \text{and} \quad \mathfrak{osp}:=\dilim{n,m}{\osp}. \]

\begin{prop} \label{prop_hamfields}
Let $(W,\omega)$ be a symplectic vector space. The map from the Lie algebra $\hl{W}$ to the Lie algebra of symplectic vector fields
\begin{displaymath}
\begin{array}{ccc}
\hl{W} & \to & \{ \eta \in \Der[\salg{W^*}] : L_\eta(\omega) = 0 \} \\
a & \mapsto & \alpha:=\Phi^{-1}(da)
\end{array}
\end{displaymath}
is a surjective homomorphism of Lie algebras with kernel $\gf$.
\end{prop}
\noproof

\begin{rem} \label{rem_hamiltonian}
We call the polynomial $a\in\hl{W}$ the \emph{Hamiltonian} corresponding to the symplectic vector field $\alpha$. By the proposition, the polynomial $a$ is uniquely determined modulo constant functions.
\end{rem}

\subsection{Lie algebra homology}

In this subsection we recall the definition of the Chevalley-Eilenberg homology of a Lie algebra.

\begin{defi} \label{def_CEhom}
Let $\mathfrak{l}$ be a Lie algebra: the underlying space of the Chevalley-Eilenberg complex of $\mathfrak{l}$ is the exterior algebra $\Lambda(\mathfrak{l})$ which is defined to be the quotient of $T(\mathfrak{l})$ by the ideal generated by the relation
\[ g \otimes h = -(-1)^{|g||h|} h \otimes g; \quad g,h \in \mathfrak{l}. \]

The differential $\delta:\Lambda^{m}(\mathfrak{l}) \to \Lambda^{m-1}(\mathfrak{l})$ is defined by the following formula:
\begin{equation} \label{eqn_CEdiff}
\delta(g_1\wedge\cdots\wedge g_m):= \sum_{1\leq i < j \leq m} (-1)^{p(g)} [g_i,g_j]\wedge g_1\wedge\cdots\wedge \hat{g_i}\wedge\cdots\wedge \hat{g_j}\wedge\cdots\wedge g_m,
\end{equation}
for $g_1,\ldots,g_m \in \mathfrak{l}$; where
\[ p(g):=|g_i|(|g_1|+\ldots+|g_{i-1}|)+|g_j|(|g_1|+\ldots+|g_{j-1}|)+|g_i||g_j|+i+j-1. \]
We will denote the Chevalley-Eilenberg complex of $\mathfrak{l}$ by $C_{\bullet}(\mathfrak{l})$. The homology of the Lie algebra $\mathfrak{l}$ is defined to be the homology of this complex. Chevalley-Eilenberg \emph{cohomology} is defined by taking the $\gf$-linear dual of the above definition.
\end{defi}

\begin{rem}
Let $\mathfrak{l}$ be a Lie algebra; $\mathfrak{l}$ acts on $\Lambda(\mathfrak{l})$ via the adjoint action. This action commutes with the Chevalley-Eilenberg differential $\delta$. As a consequence of this the space
\[ C_{\bullet}(\tglie)_{\osp} \]
of $\osp$-coinvariants of the Chevalley-Eilenberg complex of $\tglie$ forms a complex when equipped with the Chevalley-Eilenberg differential $\delta$. This is called the \emph{relative} Chevalley-Eilenberg complex of $\glie$ modulo $\osp$ (cf. \cite{fuchscohom}) and is denoted by $\cek$. The homology of this complex is called the \emph{relative} homology of $\glie$ modulo $\osp$. As before, relative \emph{cohomology} is defined by taking the $\gf$-linear dual of this definition. The stable limit of the relative Chevalley-Eilenberg complexes is given by
\[ \stlim[\bullet]=\dilim{n,m}{\cek}.\]
We denote the homology of this stable limit by $H_\bullet(\mathfrak{g},\mathfrak{osp})$ and its cohomology by $H^\bullet(\mathfrak{g},\mathfrak{osp})$.
\end{rem}

\section{Graph homology} \label{sec_grphom}

In this section we recall the definition of graph homology and formulate Kontsevich's theorem which expresses graph homology as the homology of an infinite-dimensional Lie algebra of Hamiltonian vector fields.

\subsection{Basic definitions}

This subsection recalls the basic definition of the graph complex as defined in the context of the commutative operad. Our description here will be somewhat informal, in order to avoid overburdening the reader with the precise details. A more complete description is contained in \cite{hamgraph}.

\begin{defi} \label{def_graph}
A graph $\Gamma$ is a one-dimensional cell complex. From a combinatorial perspective $\Gamma$ consists of the following data:
\begin{enumerate}
\item
A finite set, also denoted by $\Gamma$, consisting of the \emph{half-edges} of $\Gamma$.
\item \label{item_graphvert}
A partition $V(\Gamma)$ of $\Gamma$. The elements of $V(\Gamma)$ are called the \emph{vertices} of $\Gamma$.
\item \label{item_graphedge}
A partition $E(\Gamma)$ of $\Gamma$ into sets having cardinality equal to two. The elements of $E(\Gamma)$ are called the \emph{edges} of $\Gamma$.
\item \label{item_graphori}
We require an additional piece of data on $\Gamma$, namely that of an \emph{orientation}. This is given by ordering the vertices of $\Gamma$ and orienting the edges of $\Gamma$. We then mod out by the action of \emph{even} permutations on these structures, hence there are just two ways of choosing an orientation on $\Gamma$.
\end{enumerate}
\end{defi}

We say that a vertex $v \in V$ has valency $n$ if $v$ has cardinality $n$. The elements of $v$ are called the \emph{incident half-edges} of $v$. In this paper we shall assume that all our graphs have vertices of valency $\geq 3$.

There is an obvious notion of isomorphism for graphs. Two graphs $\Gamma$ and $\Gamma'$ are said to be isomorphic if there is a bijection $\Gamma\to\Gamma'$ which preserves the structures described by items (\ref{item_graphvert}--\ref{item_graphori}) of Definition \ref{def_graph}.

Given a graph $\Gamma$ and an edge $e$ of $\Gamma$, there is a way to contract the edge $e$ to produce a new graph which is denoted by $\Gamma/e$. From the topological viewpoint, this is just given by collapsing the corresponding one-cell to a point. Providing that the edge that we contract is not a loop, there is a way to make sense of the induced orientation on $\Gamma/e$.

These definitions allow us to define the following complex.

\begin{defi}
There is a complex, denoted by $\gc$, called the graph complex. The underlying space of $\gc$ is the free vector space generated by isomorphism classes of graphs modulo the relation
\[ \Gamma^\dagger = -\Gamma, \]
where $\Gamma^\dagger$ denotes the graph $\Gamma$ with the \emph{opposite} orientation.

The differential $\partial$ on $\gc$ is given by the formula:
\[ \partial(\Gamma):=\sum_{e \in E(\Gamma)} \Gamma/e \]
where $\Gamma$ is any graph and the sum is taken over all edges of $\Gamma$. Graph homology is defined to be the homology of this complex and is denoted by $H\mathcal{G}_\bullet$. Graph \emph{cohomology}, denoted by $H\mathcal{G}^\bullet$, is defined by taking the $\gf$-linear dual of the above definition.
\end{defi}

\subsection{Kontsevich's theorem}

In this subsection we formulate the super-analogue of Kontsevich's theorem linking graph homology with the homology of the Lie algebra $\glie$. A more detailed discussion of the theorem and its proof can be found in \cite{hamgraph}.

\begin{defi}
Let $c:=\{i_1,j_1\},\ldots,\{i_k,j_k\}$ be a chord diagram and assume that $i_r<j_r$ for all $r$. Let $k_1,\ldots,k_l$ be a sequence of positive integers such that $k_1+\ldots+k_l=2k$, then we can define a graph $\Gamma:=\grpchd{c}$ with half-edges $h_1,\ldots,h_{2k}$ as follows:
\begin{enumerate}
\item \label{item_grpchddummy0}
The vertices of $\Gamma$ are
\[ \{h_1,\ldots,h_{k_1}\},\{h_{k_1+1},\ldots,h_{k_1+k_2}\},\ldots,\{h_{k_1+\ldots+k_{l-1}+1},\ldots,h_{k_1+\ldots+k_l}\}. \]
\item \label{item_grpchddummy1}
The edges of $\Gamma$ are
\[ (h_{i_1},h_{j_1}),\ldots,(h_{i_k},h_{j_k}). \]
\item
The orientation on $\Gamma$ is given by ordering the vertices as in \eqref{item_grpchddummy0} and orienting the edges as in \eqref{item_grpchddummy1}.
\end{enumerate}
\end{defi}

\begin{defi} \label{def_sympinv}
Let $V$ be a vector space equipped with an \emph{even} bilinear form $\innprod$ and let
\[ c:=\{i_1,j_1\},\ldots,\{i_k,j_k\} \]
be a chord diagram, for which we assume that $i_r<j_r$ for all $r$. We define the map $\beta_c:V^{\otimes 2k}\to\gf$ by the formula
\[ \beta_c(x_1\otimes\ldots\otimes x_{2k}):= (-1)^{p_c(x)}\langle x_{i_1},x_{j_1}\rangle\langle x_{i_2},x_{j_2}\rangle\ldots\langle x_{i_k},x_{j_k}\rangle, \]
where the sign $(-1)^{p_c(x)}$ is determined by the Koszul sign rule for the permutation of the factors
\[x_1,\ldots,x_{2k} \mapsto x_{i_1},x_{j_1},x_{i_2},x_{j_2},\ldots,x_{i_k},x_{j_k}.\]
\end{defi}

We can define a map $I:\cek\to\gc$ which formally resembles Wick's formula for integration with respect to a Gaussian measure in which certain contributions to the integral are indexed by graphs:

\begin{defi} \label{def_wickmap}
Let
\[x:= (x_{11}\cdots x_{1k_1})\otimes\ldots\otimes(x_{l1}\cdots x_{lk_l})\]
represent a Chevalley-Eilenberg chain, where $x_{ij}$ is a real-valued linear function on the canonical $2n|m$-dimensional symplectic vector space $V$ (cf. \eqref{eqn_cansymp}) and $k_1+\ldots+k_l=2k$. We equip the vector space $V^*$ with the inverse form defined by diagram \eqref{fig_inverseform} and define the map
\[ I:\cek\to\gc \]
by the formula,
\[ I(x):= \sum_{c\in\chd{k}} \beta_c(x)\grpchd{c}.\]
\end{defi}

The maps $I:\cek\to\gc$ defined above give rise to a map
\begin{equation} \label{eqn_integral}
I:\stlim\to\gc
\end{equation}
on the direct limit. Kontsevich's theorem can now be formulated as follows:

\begin{theorem} \label{thm_kontthm}
The map $I:\stlim\to\gc$ is an isomorphism of complexes.
\end{theorem}
\noproof

\section{Wick's formula and the Wick rotation} \label{sec_wick}

We need to make sense of integrals of the form
\begin{equation} \label{eqn_gaussian}
\int_{\mathbb{R}^n}f(x)e^{-\frac{1}{2}\langle x,x\rangle}dx
\end{equation}
where $\innprod$ is a nondegenerate symmetric bilinear form and $f$ is a polynomial function on $\mathbb{R}^n$. If the form $\innprod$ on $\mathbb{R}^n$ is positive definite then this integral is convergent and could be evaluated by means of Wick's formula, which will be recalled below. Suppose now that $\innprod$ is \emph{not} positive definite. A linear change of variables allows one to consider only the case when the quadratic function $\langle x,x\rangle$ has the form $\sum_{i=1}^kx_i^2-\sum_{i=k+1}^nx_i^2$. We will introduce the function $g(t)$ as
\[g(t):=\int_{\mathbb{R}^n} f(x)e^{-\frac{1}{2}\left[\sum_{i=1}^k(x_i)^2+\sum_{i=k+1}^n(tx_i)^2\right]}dx.\]
Then $g(t)$ is well-defined for nonzero real $t$ and we can analytically continue $g$ for arbitrary nonzero $t\in\mathbb{C}$; thus our integral \eqref{eqn_gaussian} equals $g(i)$. For example, $\int_{-\infty}^\infty e^{\frac{1}{2}x^2}dx$ will be equal to $-i\sqrt{2\pi}$.

\begin{rem}
The formal manipulation described above is known in physics as the \emph{Wick rotation}. It is easy to check that $\int_{\mathbb{R}^n} f(x)e^{-\frac{1}{2}\langle x,x\rangle}dx$ as defined above obeys the standard rules of integration (i.e. the change of variables formula and integration by parts still hold) although the integral itself exists only formally.
\end{rem}

Now let $W$ be a real vector space of dimension $n|m$ with a nondegenerate even symmetric bilinear form $\innprod$ and denote the corresponding quadratic Hamiltonian by $\sigma:=\pi_2(\innprod)$. Given a polynomial superfunction $f\in\salg{W^*}$, we will for the most part be interested in the quantity
\begin{equation} \label{eqn_vev}
\langle f \rangle_0:=\frac{\int_W f(x,\xi)e^{-\sigma(x,\xi)}dx d\xi}{\int_W e^{-\sigma(x,\xi)}dx d\xi}
\end{equation}
which is sometimes called the \emph{vacuum expectation value} of the observable $f$ in the context of quantum field theory. Since there are no convergence issues with integrals over odd spaces, the Gaussian integrals appearing in \eqref{eqn_vev} can be defined using the Wick rotation as above by identifying the vector space $W$ with $\mathbb{R}^{n|m}$. Note that since the vacuum expectation value is a ratio of two integrals, its definition does not depend upon the choice of linear coordinates for $W$. For this reason we will hereafter omit the terms $dx$ and $d\xi$ in the formula \eqref{eqn_vev} for the expectation value. There is a purely algebraic formula, called Wick's formula, which allows one to compute  $\langle f \rangle_0$; cf. for example Lecture 1 of \cite{witten}. It suffices to assume that $f$ is a product of linear functions:

\begin{lemma} \label{lem_wicks}
Let $f=f_1\cdots f_{2k}$ be a monomial superfunction on $W$, where $f_i\in W^*$; then
\begin{equation}\label{eqn_wick}
\langle f \rangle_0=\sum_{c\in\chd{k}} \beta_c(f) = \sum_{c\in\chd{k}} (-1)^{p_c(f)}\langle f_{i_1}, f_{j_1} \rangle^{-1} \ldots \langle f_{i_k}, f_{j_k} \rangle^{-1};
\end{equation}
where the sum is taken over all partitions $c:=\{i_1,j_1\},\ldots,\{i_k, j_k\}$ of $\{1,\ldots, 2k\}$, also known as chord diagrams. The sign $(-1)^{p_c(f)}$ is determined by the Koszul sign rule (cf. Definition \ref{def_sympinv}).
\end{lemma}
\noproof

\begin{rem}
It follows that $\langle f \rangle_0$ is always real and could in fact be defined over an arbitrary field of characteristic zero.
\end{rem}

We now provide the appropriate analogue of Stokes' theorem.

\begin{lemma} \label{lem_stokes}
Let $\innprod$ be a nondegenerate even symmetric bilinear form on $\gf^n$ and let $q\in\forms[n-1]{\gf^n}$ be any polynomial $(n-1)$-form, then
\[\int_{\gf^n} d\left(q e^{-\sigma}\right) = 0.\]
\end{lemma}

\begin{proof}
It suffices to consider the case $n=1$. We may assume that the function $q(x)=x^r$ is a monomial in one variable. First let us suppose that our bilinear form is \emph{positive} definite, then
\[ \int_{-\infty}^{\infty} \pder{x}\left[x^r e^{-\frac{1}{2}x^2}\right]dx = (1-(-1)^r)\lim_{a\to\infty} a^r e^{-\frac{1}{2}a^2}=0.\]
The case when $\innprod$ is \emph{negative} definite can be reduced to the first case by making use of the Wick rotation.
\end{proof}

\section{The geometry of the BV-formalism} \label{sec_BVform}

In this section we describe the basic geometry underlying the BV-formalism. This description is not new; it is taken from \cite{schwarz}, but also appears in the work of several other authors, cf. \cite{Kh} and \cite{KS}. Indeed, our situation is considerably simpler, since we only work with spaces with a flat geometry. The only significant difference is that our spaces are no longer compact and that we use Gaussian measures defined from a nondegenerate bilinear form which may not be positive definite.

\begin{defi}
A \emph{linear $P$-manifold} is a symplectic vector space $(W,\omega)$ whose symplectic form $\omega\in\forms[2]{W}$ is \emph{odd}.
\end{defi}

Since the bilinear form on $W$ corresponding to $\omega$ must be nondegenerate, it follows that any linear $P$-manifold must have dimension $n|n$ for some integer $n$.

An important fact which we will use is that for every positive integer $n$, there is a canonical linear $P$-manifold $U$ of dimension $n|n$ to which any other $n|n$-dimensional linear $P$-manifold is isomorphic. There are canonical coordinates
\[ \underbrace{x_1,\ldots,x_n}_{\text{even}};\underbrace{\xi_1,\ldots,\xi_n}_{\text{odd}} \]
on $U$ which form a basis for $U^*$ and the symplectic form $\omega$ is given by
\begin{equation} \label{eqn_canPman}
\omega:=\sum_{i=1}^n dx_i d\xi_i.
\end{equation}

\subsection{Divergence}

We now give the definition of (super)divergence which generalises the definition of (super)trace.

\begin{defi}
Let $W$ be a finite-dimensional vector space, then we define a map
\[ \nabla:\Der[\salg{W^*}]\to\salg{W^*} \]
as follows: let $y_1,\ldots,y_k$ be a homogeneous basis for $W^*$ and $\eta:\salg{W^*}\to\salg{W^*}$ be any vector field, then
\[\nabla\eta:=\sum_{i=1}^k (-1)^{y_i+y_i\eta}\pder{y_i}[\eta(y_i)].\]
\end{defi}

\begin{rem}
The definition of $\nabla$ is independent of the choice of linear coordinates as it is given by the composition of the maps in the following diagram:
\[\Der[\salg{W^*}]=\Hom(W^*,\salg{W^*})=\salg{W^*}\otimes W\overset{\tau}{\to} W\otimes\salg{W^*}\to\salg{W^*},\]
where the last map is the action of $W=W^{**}$ on $\salg{W^*}$ given by considering $W$ as the space of constant vector fields.
\end{rem}

The divergence operator $\nabla$ can easily be shown to satisfy the following identities:

\begin{lemma} \label{lem_dividentities}
Let $W$ be a finite-dimensional vector space; $\eta,\gamma:\salg{W^*}\to\salg{W^*}$ be vector fields and $f\in\salg{W^*}$ be a function:
\begin{enumerate}
\item \label{item_dividentitiescomm}$\nabla[\eta,\gamma]=\eta(\nabla\gamma)-(-1)^{\eta\gamma}\gamma(\nabla\eta)$,
\item \label{item_dividentitiesfmult} $\nabla(f\cdot\eta)=f\cdot\nabla\eta+(-1)^{f\eta}\eta(f)$.
\end{enumerate}
\end{lemma}
\noproof

\subsection{Anti-bracket}

Recall from Lemma \ref{lem_fieldsforms} that given a symplectic vector space $(W,\omega)$, there corresponds to any polynomial function $a\in\salg{W^*}$, a symplectic vector field $\alpha:=\Phi^{-1}(da)$. We now use this fact to define the anti-bracket.

\begin{defi}
Let $(W,\omega)$ be a linear $P$-manifold. We define a bracket
\[ \{-,-\}:\salg{W^*}\otimes\salg{W^*}\to\salg{W^*} \]
of \emph{odd} degree by the formula
\[ \{a,b\}:=L_\alpha(b).\]
\end{defi}

\begin{lemma} \label{lem_lieids}
The bracket $\{-,-\}$ on $\hl{W}:=\salg{W^*}$ is an \emph{odd} Poisson bracket, that is to say that:
\begin{enumerate}
\item
The bracket $[-,-]:\Pi\hl{W}\otimes\Pi\hl{W}\to\Pi\hl{W}$ given by the formula
\[ \Pi\circ[-,-] = \{-,-\}\circ(\Pi\otimes\Pi) \]
is a Lie bracket.
\item \label{item_lieids2}
For any $a,b,c\in\hl{W}$;
\[\{a,bc\}=\{a,b\}c+(-1)^{(a+1)b}b\{a,c\}.\]
\end{enumerate}
\end{lemma}

\begin{proof}
Firstly, the bracket $\{-,-\}$ is symmetric and hence the bracket $[-,-]$ is anti-symmetric; this is a simple consequence of identities \eqref{item_operatoridsone} and \eqref{item_operatoridsfour} of Lemma \ref{lem_operatorids}.

The Jacobi identity for $[-,-]$ follows from the following identity:
\begin{equation} \label{eqn_bracketjacobi}
\begin{split}
d\{a,b\} &= dL_\alpha(b) = (-1)^{a+1}L_\alpha d(b) = (-1)^{a+1}L_\alpha i_\beta(\omega), \\
&=(-1)^{a+1}[L_\alpha,i_\beta](\omega) = (-1)^{a+1}i_{[\alpha,\beta]}(\omega). \\
\end{split}
\end{equation}

Finally, the Leibniz identity follows from the identity
\begin{equation} \label{eqn_bracketleibniz}
d(bc) = (-1)^{c(b+1)}cd(b) + (-1)^{b}bd(c) = \Phi((-1)^{c(b+1)}c\beta + (-1)^{b}b\gamma).
\end{equation}
\end{proof}

\begin{rem} \label{rem_oddhamfields}
As in section \ref{sec_sympgeom}, we define $\gl{W}$ to be the Lie subalgebra of $\hl{W}$ consisting of polynomials of order $\geq 2$. The corresponding analogue of Proposition \ref{prop_hamfields} holds for $\gl{W}$; that is to say that the map
\[[\Phi^{-1}d\Pi]:\Pi\gl{W}\to\{\eta\in\Der_+[\salg{W^*}]:L_\eta(\omega)=0\}\]
is an isomorphism of Lie algebras.
\end{rem}

\subsection{Odd Laplacian} \label{sec_BVlap}

Now we define the (odd) Laplacian.

\begin{defi}
Let $(W,\omega)$ be a linear $P$-manifold. We define the Laplacian
\[ \Delta:\salg{W^*}\to\salg{W^*} \]
by the formula
\[ \Delta(a):=\frac{1}{2}\nabla(\alpha).\]
\end{defi}

\begin{lemma} \label{lem_BVidentity}
The Laplacian $\Delta:\salg{W^*}\to\salg{W^*}$ satisfies the following identities:
\begin{enumerate}
\item \label{item_BVidentity1}
For all $a,b\in\salg{W^*}$;
\[\Delta(ab)=\Delta(a)b + (-1)^{a}a\Delta(b) + \{a,b\}.\]
\item \label{item_BVidentity2}
For all $a,b\in\salg{W^*}$;
\[\Delta\{a,b\} + \{\Delta a,b\} +(-1)^{a}\{a,\Delta b\} = 0.\]
\item \label{item_BVidentity3}
Suppose that $x_1,\ldots,x_n;\xi_1,\ldots,\xi_n$ is a system of coordinates on $W$ such that $\omega$ takes the canonical form \eqref{eqn_canPman}, then
\begin{equation} \label{eqn_canlap}
\Delta=\sum_{i=1}^n\pder{x_i}\pder{\xi_i}.
\end{equation}
\end{enumerate}
\end{lemma}

\begin{rem}
Note that it follows immediately from \eqref{eqn_canlap} that $\Delta^2=0$.
\end{rem}

\begin{proof}
Part \eqref{item_BVidentity1} follows directly from Lemma \ref{lem_dividentities} \eqref{item_dividentitiesfmult} and equation \eqref{eqn_bracketleibniz}. Part \eqref{item_BVidentity2} is a direct consequence of Lemma \ref{lem_dividentities} \eqref{item_dividentitiescomm} and equation \eqref{eqn_bracketjacobi}. The proof of part \eqref{item_BVidentity3} is a routine calculation.
\end{proof}

\subsection{Lagrangian subspaces and Stokes' theorem}

The integration of functions over Lagrangian submanifolds of a $P$-manifold lies at the heart of the BV-formalism. It was a fundamental insight of Schwarz \cite{schwarz} that if the $P$-manifold comes equipped with a compatible volume form ($SP$-manifold), then any Lagrangian submanifold may also be endowed with a volume form which is canonically determined up to a sign. Moreover, the integral of a $\Delta$-closed superfunction over any two homologous submanifolds is always the same; this is the gauge-independence property of the BV-framework.

In this subsection we discuss linear analogues of some of these results. Our integrals are defined over a noncompact space supplied with a Gaussian measure which is determined by a (not necessarily positive definite) bilinear form. Because of this we cannot formally quote the relevant statements from Schwarz's paper and therefore provide an independent treatment.

\begin{defi}
Let $(W,\omega)$ be a linear $P$-manifold. A subspace $X\subset W$ is called \emph{isotropic} if the restriction of the symplectic form $\omega\in\forms[2]{W}$ to $X$ is equal to zero. A subspace $L\subset W$ is called a \emph{Lagrangian subspace} if it is a \emph{maximally isotropic} subspace of $W$; that is to say that it is an isotropic subspace of maximal dimension.
\end{defi}

Let $U$ be the canonical linear $P$-manifold of dimension $n|n$ (cf. \eqref{eqn_canPman}). There is a canonical Lagrangian subspace $L_{k|n-k}\subset U$ of dimension $k|n-k$ which is defined as the locus of the equations
\begin{equation} \label{eqn_lagsubspace}
\xi_1=\ldots=\xi_k=0,\quad x_{k+1}=\ldots=x_n=0.
\end{equation}

We now formulate some basic facts about Lagrangian subspaces.

\begin{lemma} \label{lem_lagsubspace}
Let $(W,\omega)$ be a linear $P$-manifold of dimension $n|n$ and let $L\subset W$ be a Lagrangian subspace:
\begin{enumerate}
\item \label{item_lagsubspace1}
There exist canonical coordinates $x_1,\ldots,x_n;\xi_1,\ldots,\xi_n$ on $W$ such that $\omega$ takes the canonical form \eqref{eqn_canPman} and $L$ is the locus of the equations \eqref{eqn_lagsubspace} for some $k \leq n$.
\item \label{item_lagsubspace2}
There exists a Lagrangian subspace $L'\subset W$ such that
\[ W=L\oplus L'.\]
\item \label{item_lagsubspace3}
For any such Lagrangian subspace $L'$, the map
\begin{displaymath}
\begin{array}{ccc}
L' & \to & L^* \\
x & \mapsto & [y\mapsto\langle x,y\rangle] \\
\end{array}
\end{displaymath}
is nondegenerate.
\item \label{item_lagsubspace4}
Any Lagrangian subspace can be represented as the graph of a gradient; more precisely, suppose that $x_1,\ldots,x_n;\xi_1,\ldots,\xi_n$ is a system of coordinates on $W$ for which $\omega$ takes the canonical form \eqref{eqn_canPman} and such that the coordinates
\[x_1,\ldots,x_k;\xi_{k+1},\ldots,\xi_n\]
are linearly independent on $L$, then there exists an odd quadratic function (known as a generating function) $\varphi\in S^2(W^*)$ such that
\[\pder{\xi_1}(\varphi)=\ldots=\pder{\xi_k}(\varphi)=0=\pder{x_{k+1}}(\varphi)=\ldots=\pder{x_n}(\varphi)\]
and such that $L$ is the locus of the equations
\[\xi_1=-\pder{x_1}(\varphi),\ldots,\xi_k=-\pder{x_k}(\varphi),\quad x_{k+1}=\pder{\xi_{k+1}}(\varphi),\ldots, x_n=\pder{\xi_n}(\varphi).\]
\end{enumerate}
\end{lemma}

\begin{rem}
It follows that the dimension of any Lagrangian subspace $L$ is then $k|n-k$ for some $k\leq n$. Note that any subspace given by the graph of a gradient is Lagrangian, as a routine calculation using the equality of mixed partial derivatives shows. This implies that the Lagrangian subspace $L'$ in part \eqref{item_lagsubspace2} is not canonically determined, the different choices for $L'$ being parameterised by some chart of the Lagrangian Grassmannian $\mathcal{L}(W,\omega)$.
\end{rem}

\begin{proof}
The proof of part \eqref{item_lagsubspace1} is an exercise in linear algebra which is omitted. Part \eqref{item_lagsubspace2} follows from part \eqref{item_lagsubspace1}; if $L$ is the locus of the equations \eqref{eqn_lagsubspace} then we can set $L'$ to be the locus of the equations
\[x_1=\ldots=x_k=0,\quad \xi_{k+1}=\ldots=\xi_n=0.\]
The proof of part \eqref{item_lagsubspace3} is tautological. The proof of part \eqref{item_lagsubspace4} follows the proof of the corresponding classical result mutatis mutandis, cf. \cite{mcduff}.
\end{proof}

The following theorem describes a kind of Poincar\'e duality which allows us to identify the geometry underlying the BV-formalism with the usual exterior calculus.

\begin{theorem} \label{thm_poincareduality}
Let $U$ be the canonical linear $P$-manifold of dimension $n|n$ (cf. \eqref{eqn_canPman}) and let $M:=U_0$ be the body (even component) of $U$. There is a duality isomorphism (inverse Fourier transform)
\[ \mathcal{D}:\csalg{U^*}\to\cforms{M} \]
given as follows: let $g:=f\xi_{i_1}\ldots\xi_{i_l}$ where $f\in\csalg{M^*}$ and $\xi_{i_r}\in U_1^*$, then
\[\mathcal{D}(g):= i_{\pder{x_{i_1}}}\circ\ldots\circ i_{\pder{x_{i_l}}}[f\cdot dx_1\ldots dx_n].\]
\end{theorem}

\begin{proof}
It follows from Lemma \ref{lem_operatorids} \eqref{item_operatoridsfour} that $\mathcal{D}$ is well-defined. The statement that $\mathcal{D}$ is an isomorphism is simply a reformulation of the usual statement of Poincar\'e duality which holds for the exterior algebra. It remains to show that $\mathcal{D}$ identifies the Laplacian $\Delta$ with the exterior derivative $d$:
\begin{displaymath}
\begin{split}
\Delta(g) &= \sum_{r=1}^l (-1)^{r-1}\pder{x_{i_r}}(f)\xi_{i_1}\ldots\widehat{\xi_{i_r}}\ldots\xi_{i_l}, \\
d\mathcal{D}(g) &= \sum_{r=1}^l(-1)^{r-1} i_{\pder{x_{i_1}}}\circ\ldots \widehat{i_{\pder{x_{i_r}}}}\ldots\circ i_{\pder{x_{i_l}}}\circ L_{\pder{x_{i_r}}}[f\cdot dx_1\ldots dx_n]. \\
\end{split}
\end{displaymath}
The first identity follows from Lemma \ref{lem_BVidentity} \eqref{item_BVidentity3}. The second identity follows from parts \eqref{item_operatoridsone} and \eqref{item_operatoridstwo} of Lemma \ref{lem_operatorids} and the fact that $f\cdot dx_1\ldots dx_n$ is a top form and therefore closed. The above identities imply that $\mathcal{D}\Delta(g)=d\mathcal{D}(g)$.
\end{proof}

We have the following proposition whose proof is a routine check.

\begin{prop} \label{prop_intduality}
Let $U$ be the canonical linear $P$-manifold of dimension $n|n$ and let $M:=U_0$ be the $n$-dimensional body of $U$. Let $L_{k|n-k}$ be the canonical Lagrangian subspace of $U$ of dimension $k|n-k$ (cf. \eqref{eqn_lagsubspace}) and denote the $k$-dimensional body of $L_{k|n-k}$ by $M_k$. Let $f\in\salg{U^*}$ be a polynomial superfunction on $U$ and let $\sigma\in S^2(U^*)$ be a quadratic even Hamiltonian whose restriction to $L_{k|n-k}$ is nondegenerate, then
\begin{equation} \label{eqn_intduality}
\int_{L_{k|n-k}} f e^{-\sigma}\, dx_1\ldots dx_k d\xi_{k+1}\ldots d\xi_n = (-1)^{k(n-k)}\int_{M_k} \mathcal{D}(f e^{-\sigma}).
\end{equation}
\end{prop}

\begin{rem}
Here the operator $\int_{M_k}$ is defined to be zero on $i$-forms for $i<k$.
\end{rem}
\noproof

\begin{cor} \label{cor_BVstokes}
Let $f\in\salg{U^*}$ be any polynomial on $U$ and let $\sigma\in S^2(U^*)$ be a quadratic even Hamiltonian whose restriction to $L_{k|n-k}$ is nondegenerate, then
\[ \int_{L_{k|n-k}} \Delta(f e^{-\sigma})\,dx_1\ldots dx_k d\xi_{k+1}\ldots d\xi_n=0. \]
\end{cor}

\begin{proof}
This is a direct consequence of Lemma \ref{lem_stokes}, Theorem \ref{thm_poincareduality} and Proposition \ref{prop_intduality}.
\end{proof}

\section{The dual construction} \label{sec_dual}

In this section we describe the `dual construction' introduced by Kontsevich in \cite{kontfd} which produces cohomology classes in the graph complex from a contractible differential graded commutative Frobenius algebra. We provide a formulation of this construction in terms of the BV-formalism and use this framework to show the independence of the obtained cohomology class on the choice of a gauge-condition.

\subsection{Preliminary combinatorics}

Given an $\ai$-structure on a vector space and an associative algebra, one can form an $\ai$-structure on the tensor product in a natural way. This construction can be extended to all coderivations (not necessarily those determining an $\ai$-structure) to define the following map.

\begin{defi}
Let $W$ be any vector space and $A$ be an associative algebra. Let $m_n:A^{\otimes n}\to A$ be the map defined by the formula
\[m_n(a_1,\ldots,a_n):=a_1\cdots a_n,\quad n\geq 1 .\]
Using the correspondence between coderivations and multilinear maps provided by \eqref{eqn_correspondence} we can define a map $\Psi_A:\Coder_+[T(W)]\to\Coder_+[T(A\otimes W)]$ which is given by the following commutative diagram:
\[\xymatrix{A^{\otimes n}\otimes W^{\otimes n} \ar@{=}[r] \ar[d]_{m_n\otimes\zeta_n} & (A \otimes W)^{\otimes n} \ar[ld]^{\Psi_A(\zeta_n)} \\ A\otimes W}\]
for any map $\zeta_n:W^{\otimes n}\to W$.
\end{defi}

\begin{defi}
Let $\zeta_n:W^{\otimes n}\to W$ and $\gamma_m:W^{\otimes m}\to W$ be maps, we define a map
\[ \zeta_n\circ_i\gamma_m:W^{\otimes n+m-1}\to W \]
for $1\leq i\leq n$ by the formula
\[\zeta_n\circ_i\gamma_m(x_1,\ldots,x_{n+m-1}):= (-1)^{(x_1+\ldots+x_{i-1})\gamma_m}\zeta_n(x_1,\ldots,x_{i-1},\gamma_m(x_i,\ldots,x_{i+m-1}),x_{i+m},\ldots,x_{n+m-1}).\]
\end{defi}

\begin{rem}
In this notation the formula for the commutator bracket can be written as follows:
\begin{equation} \label{eqn_commbracket}
[\zeta_n,\gamma_m]:= \sum_{i=1}^n \zeta_n\circ_i\gamma_m -(-1)^{\zeta_n\gamma_m}\sum_{j=1}^m\gamma_m\circ_j\zeta_n.
\end{equation}
\end{rem}

\begin{lemma} \label{lem_mapliealg}
Let $W$ be a vector space and let $A$ be an associative algebra, then the map
\[ \Psi_A:\Coder_+[T(W)]\to\Coder_+[T(A\otimes W)] \]
is a map of Lie algebras.
\end{lemma}

\begin{proof}
This follows from the fact that $m_n\circ_i m_k=m_{n+k-1}$ for all $i$, giving us the identity
\begin{displaymath}
\begin{split}
\Psi_A(\zeta_n)\circ_i\Psi_A(\gamma_k) &= (m_n\otimes\zeta_n)\circ_i(m_k\otimes\gamma_k), \\
&= (m_n\circ_i m_k)\otimes(\zeta_n\circ_i\gamma_k), \\
&= \Psi_A(\zeta_n\circ_i\gamma_k). \\
\end{split}
\end{displaymath}
for all maps $\zeta_n:W^{\otimes n}\to W$ and $\gamma_k:W^{\otimes k}\to W$. Applying this identity to the formula for the commutator bracket given by equation \eqref{eqn_commbracket} yields the desired result.
\end{proof}

Let $A$ be a differential graded associative algebra and let $W$ be any vector space. We can define a differential
\[ \tilde{d}:A\otimes W\to A\otimes W \]
by the formula $\tilde{d}:= d\otimes\id$. We have the following lemma.

\begin{lemma} \label{lem_mapleibniz}
For any map $\zeta_n:W^{\otimes n}\to W$,
\[[\tilde{d},\Psi_A(\zeta_n)]=0.\]
\end{lemma}

\begin{proof}
Since $d$ is a derivation we have
\[[d,m_n]= d\circ_1 m_n - \sum_{i=1}^n m_n\circ_i d = 0.\]
Now
\begin{displaymath}
\begin{split}
[\tilde{d},\Psi_A(\zeta_n)] &= \tilde{d}\circ_1 \Psi_A(\zeta_n) - (-1)^{\zeta_n}\sum_{i=1}^n \Psi_A(\zeta_n)\circ_i \tilde{d} , \\
&= (d\otimes\id)\circ_1 (m_n\otimes\zeta_n) - (-1)^{\zeta_n}\sum_{i=1}^n (m_n\otimes\zeta_n)\circ_i (d\otimes\id) , \\
&= [d,m_n]\otimes\zeta_n =0.\\
\end{split}
\end{displaymath}
\end{proof}

Now we need to formulate the commutative analogue of the above. Let $W$ be any vector space and equip the symmetric algebra $S(W)$ with the canonical cocommutative diagonal. By using the correspondence \eqref{eqn_correspondence} between coderivations and multilinear maps, we can define a surjective map of Lie algebras
\[ \xymatrix{\Coder[T(W)] \ar@{->>}[rrr]^{P_W} &&& \Coder[S(W)] \\ \prod_{n=0}^\infty \Hom(W^{\otimes n},W) \ar@{=}[u] \ar[rrr]^{\zeta_n\mapsto\zeta_n\circ i_n} &&& \prod_{n=0}^\infty \Hom(S^n(W),W) \ar@{=}[u]} \]

\begin{rem}
Any $\ai$-algebra could be regarded as a nilpotent element in $\Coder[T(W)]$. Its image under the map $P_W$ corresponds to the associated commutator $L_\infty$-algebra.
\end{rem}

\begin{defi}
Let $A$ be a \emph{commutative} algebra and define $m_n:S^n(A)\to A$ by the formula
\[ m_n(a_1,\ldots,a_n):= a_1\cdots a_n.\]
For any vector space $W$ we define the map $\Psi_A:\Coder_+[S(W)]\to\Coder_+[S(A\otimes W)]$ by the following commutative diagram:
\[\xymatrix{S^n(A)\otimes S^n(W) \ar[d]_{m_n\otimes\zeta_n} & S^n(A \otimes W) \ar@{->>}[l] \ar[ld]^{\Psi_A(\zeta_n)} \\ A\otimes W}\]
\end{defi}

\begin{lemma} \label{lem_comliemap}
The maps defined above fit into the following commutative diagram:
\[ \xymatrix{ \Coder_+[S(W)] \ar[rr]^{\Psi_A} && \Coder_+[S(A\otimes W)] \\ \Coder_+[T(W)] \ar[u]^{P_W} \ar[rr]^{\Psi_A} && \Coder_+[T(A\otimes W)] \ar[u]^{P_{A\otimes W}}} \]
\end{lemma}

\begin{proof}
This follows as a straightforward consequence of the commutativity of the algebra $A$.
\end{proof}

\begin{cor} \label{cor_leibnizvanishing}
Let $W$ be a vector space and $A$ be a commutative differential graded algebra:
\begin{enumerate}
\item
The map $\Psi_A:\Coder_+[S(W)]\to\Coder_+[S(A\otimes W)]$ is a map of Lie algebras,
\item \label{item_leibnizvanishing2}
For all $\zeta_n:S^n(W)\to W$,
\[ [\tilde{d},\Psi_A(\zeta_n)]=0.\]
\end{enumerate}
\end{cor}

\begin{proof}
This follows as a direct consequence of Lemmas \ref{lem_mapliealg}, \ref{lem_mapleibniz} and \ref{lem_comliemap} and the fact that the map $P_W:\Coder_+[T(W)]\to\Coder_+[S(W)]$ is surjective.
\end{proof}

Since diagram \eqref{eqn_coderder} gives us an identification between coderivations and derivations, one can consider the map $\Psi_A$ as a map between Lie algebras of formal vector fields. Since this map preserves the order of these vector fields, it restricts to a map on polynomial vector fields. These facts are expressed by the following commutative diagram:
\[ \xymatrix{ \Coder_+[\salg{A\otimes W}] \ar[rr]^{\zeta\mapsto\zeta^{\vee}}_{\cong} && \Der_+[\csalg{A^*\otimes W^*}] && \Der_+[\salg{A^*\otimes W^*}] \ar@{_{(}->}[ll]_{\hat{\eta}\mapsfrom\eta} \\ \Coder_+[\salg{W}] \ar[rr]^{\zeta\mapsto\zeta^{\vee}}_{\cong} \ar[u]^{\Psi_A} && \Der_+[\csalg{W^*}] \ar[u]^{\Psi_A} && \Der_+[\salg{W^*}] \ar@{_{(}->}[ll]_{\hat{\eta}\mapsfrom\eta} \ar[u]^{\Psi_A} } \]

Suppose that $W$ is a vector space with a nondegenerate skew-symmetric bilinear form $\langle-,-\rangle_W$ (i.e. a symplectic vector space) and $A$ is a vector space with a nondegenerate symmetric bilinear form $\langle-,-\rangle_A$, then we can define a nondegenerate skew-symmetric bilinear form $\langle-,-\rangle_{A\otimes W}$ on $A\otimes W$ by tensoring these bilinear forms together:
\[ \langle a_1\otimes w_1, a_2\otimes w_2\rangle_{A\otimes W}:= (-1)^{w_1 a_2}\langle a_1,a_2\rangle_A\langle w_1,w_2\rangle_W. \]

Now suppose that $A$ is a commutative Frobenius algebra and that the bilinear forms $\innprod_W$ and $\innprod_A$ are even and odd respectively. We can extend the definition of $\Psi_A$ to the corresponding Lie algebras of Hamiltonians. Let $\tilde{m}_n \in (A^*)^{\otimes n+1}$ be the tensor given by the formula
\begin{equation} \label{eqn_itmult}
\tilde{m}_n(a_1,\ldots,a_{n+1}):=\langle m_n(a_1,\ldots,a_n),a_{n+1}\rangle_A.
\end{equation}
It follows from condition \eqref{eqn_invinnprod} and the fact that $A$ is commutative that $\tilde{m}_n$ is an $S_{n+1}$-invariant tensor.

\begin{lemma}
There exists a unique map $\Psi_A:\gl{W}\to\gl{A\otimes W}$ making the following two diagrams commute:
\[\xymatrix{ S^{n+1}(W^*) \ar[rrr]^{\Psi_A} &&& S^{n+1}(A^*\otimes W^*) \\ (W^*)^{\otimes n+1} \ar@{->>}[u] \ar[rr]^-{x\mapsto \tilde{m}_n\otimes x} && (A^*)^{\otimes n+1}\otimes (W^*)^{\otimes n+1} \ar@{=}[r] & (A^*\otimes W^*)^{\otimes n+1} \ar@{->>}[u] }\]
\[\xymatrix{ \gl{A\otimes W} \ar[rrr]^{[\Phi_{A\otimes W}^{-1}\circ d]} &&& \Der_+[\salg{A^*\otimes W^*}] \\ \gl{W} \ar[rrr]^{[\Phi_W^{-1}\circ d]} \ar[u]^{\Psi_A} &&& \Der_+[\salg{W^*}] \ar[u]^{\Psi_A} } \]
where the map $\Phi$ between vector fields and 1-forms was defined in Lemma \ref{lem_fieldsforms}.

Furthermore, the map
\[[\Pi\Psi_A]:\gl{W}\to\Pi\gl{A\otimes W}\]
is a map of Lie algebras.
\end{lemma}

\begin{rem}
Note that the map $\Psi_A$ is not a map of commutative algebras.
\end{rem}

\begin{proof}
Consider the first diagram. The map on the bottom row is well-defined and lifts uniquely to a map $\Psi_A$ on the top row since the tensor $\tilde{m}_n$ is $S_{n+1}$-invariant. The commutativity of the second diagram follows from a relatively straightforward calculation, which again uses the fact that the tensor $\tilde{m}_n$ is $S_{n+1}$-invariant. It now follows from Proposition \ref{prop_hamfields} (cf. also Remark \ref{rem_oddhamfields}) and Corollary \ref{cor_leibnizvanishing} that $\Pi\Psi_A$ is a map of Lie algebras.
\end{proof}

\begin{prop} \label{prop_divsupertrace}
Let $W$ be a vector space and let $\zeta_n:S^n(W)\to W$ be a map. We have the following formula for the divergence $[i_{n-1}\nabla(\zeta_n^{\vee})]:S^{n-1}(W)\to\gf$ of this vector field:
\begin{equation} \label{eqn_divsupertrace}
a_1,\ldots,a_{n-1} \mapsto \tr[x\mapsto\zeta_n(a_1,\ldots,a_{n-1},x)].
\end{equation}
\end{prop}

\begin{proof}
We may assume that we can write the map $\zeta_n$ as
\[ \zeta_n:= y\otimes\alpha\in [S^n(W)]^*\otimes W. \]
The element $\alpha\in W$ may be considered as a constant vector field on $\salg{W^*}$; from this we get
\[ \nabla(\zeta_n^{\vee}) = (-1)^{y\alpha}\alpha\pi_n(y)=(-1)^{y\alpha}\frac{1}{n}\pi_{n-1}\alpha(y).\]
Consider the map
\[ T:=[\id^{\otimes n-1}\otimes\tr]:(W^*)^{\otimes n}\otimes W \to (W^*)^{\otimes n-1}. \]
By a straightforward calculation one can prove that
\[(-1)^{y\alpha}\alpha(y)=\sum_{i=1}^n z_{n-1}^{i-1}T[(z_n^{-i}\cdot y)\otimes\alpha],\]
where $z_n$ is the cyclic permutation $(1 \, 2 \ldots n)$. Since the tensor $y$ is symmetric, it follows that
\[ \nabla(\zeta_n^{\vee})= \pi_{n-1}T(y\otimes\alpha)=\pi_{n-1}T(\zeta_n). \]
The proposition now follows from the fact that the tensor $T(\zeta_n)$ is symmetric and represents the map \eqref{eqn_divsupertrace}.
\end{proof}

We now give a proof of the following crucial theorem. This theorem is peculiar to \emph{commutative} geometry, i.e. it only holds in the context of graphs decorated by the commutative operad and even then, only with the aforementioned notion of orientation. Recall that the commutative graph complex does not allow graphs with loops; the reason for this is that such graphs admit an orientation reversing automorphism given by permuting the half-edges of the corresponding loop, which thus means that they must vanish in the graph complex. The proof of the vanishing of a divergence like quantity in the following theorem exploits the same kind of symmetry.

\begin{theorem} \label{thm_vanishing}
Let $W$ be a vector space with a nondegenerate even skew-symmetric bilinear form $\innprod_W$ and let $A$ be a commutative Frobenius algebra whose nondegenerate bilinear form $\innprod_A$ is symmetric and odd. Let $\zeta_n:S^n(W)\to W$ be any map and set $\gamma_n:=\Psi_A(\zeta_n)$, then the following identity holds:
\[ \nabla\left(\gamma_n^{\vee}\right) = 0.\]
\end{theorem}

\begin{proof}
Let $\tilde{\gamma}_n \in [A\otimes W]^{*\otimes n+1}$ be the map defined by the formula
\[ \tilde{\gamma}_n(z_1,\ldots,z_{n+1}):=\langle\gamma_n(z_1,\ldots,z_n),z_{n+1}\rangle_{A\otimes W}. \]
Consider the map
\[ T:=\left[\id^{\otimes n-1}\otimes\innprod_{A\otimes W}^{-1}\right]:[A\otimes W]^{*\otimes n+1}\to[A\otimes W]^{*\otimes n-1}. \]
A straightforward calculation yields the identity
\[ [T(\tilde{\gamma}_n)](z_1,\ldots,z_{n-1}) = (-1)^{\gamma_n}\tr[x\mapsto\gamma_n(z_1,\ldots,z_{n-1},x)]; \]
hence it follows from Proposition \ref{prop_divsupertrace} that $T(\tilde{\gamma}_n)=(-1)^{\gamma_n}i_{n-1}\nabla\left(\gamma_n^{\vee}\right)$.

Now define the map $\tau:[A^*\otimes W^*]^{\otimes 2}\to[A^*\otimes W^*]^{\otimes 2}$ to be the map given by permuting the factors of $A$ and leaving the factors of $W$ fixed:
\[ \tau([f_1\otimes g_1]\otimes[f_2\otimes g_2]):= (-1)^{f_2g_1+f_2f_1+f_1g_1}[f_2\otimes g_1]\otimes[f_1\otimes g_2], \]
where $f_1,f_2\in A^*$ and $g_1,g_2 \in W^*$.

Since the tensor \eqref{eqn_itmult} is \emph{symmetric}, it follows that
\begin{equation} \label{eqn_vanishingdummy1}
\tilde{\gamma}_n = \left[\id^{\otimes n-1}\otimes\tau\right](\tilde{\gamma}_n).
\end{equation}

On the other hand, a straightforward calculation gives the following identity for the inverse form:
\[ \innprod_{A\otimes W}^{-1}=\innprod_A^{-1}\otimes\innprod_W^{-1}. \]
Note that since the bilinear form $\innprod_A$ is \emph{odd}, it follows that the corresponding inverse form $\innprod_A^{-1}$ is \emph{anti-symmetric}, due to the presence of signs coming from the Koszul sign rule in diagram \eqref{fig_inverseform}. From this fact it follows that
\begin{equation} \label{eqn_vanishingdummy2}
\innprod_{A\otimes W}^{-1} = -\left[\innprod_{A\otimes W}^{-1}\circ\tau\right].
\end{equation}

Combining \eqref{eqn_vanishingdummy1} and \eqref{eqn_vanishingdummy2} we see that
\[ T(\tilde{\gamma}_n) =  T\circ \left[\id^{\otimes n-1}\otimes\tau\right](\tilde{\gamma}_n) = \left[\id^{\otimes n-1}\otimes\left[\innprod_{A\otimes W}^{-1}\circ\tau\right]\right](\tilde{\gamma}_n) = -T(\tilde{\gamma}_n) \]
and therefore we must have $T(\tilde{\gamma}_n)=0$.
\end{proof}

\subsection{The approach through the BV-formalism} \label{sec_dualBVform}

Throughout the rest of this subsection $A$ will denote a commutative differential graded Frobenius algebra. Furthermore, we will assume that the nondegenerate bilinear form $\innprod_A$ is \emph{odd} and \emph{symmetric} and that the differential $d$ is \emph{contractible}. We can define a new (degenerate) anti-symmetric bilinear form $\innprod_d$ on $A$ by the formula
\begin{equation} \label{eqn_degform}
\langle a_1,a_2 \rangle_d:= (-1)^{a_1}\langle a_1,d(a_2) \rangle_A.
\end{equation}

Let $W$ be a symplectic vector space with an \emph{even} symplectic form. Recall that $A\otimes W$ will become a symplectic vector space with an \emph{odd} symplectic form and that the induced differential on $A\otimes W$ is denoted by $\tilde{d}:=d\otimes\id$. It follows from condition \eqref{eqn_dinv} that $\tilde{d}^{\vee}$ is a linear symplectic vector field and hence we can consider its corresponding even quadratic Hamiltonian $\tilde{\sigma}\in S^2(A^*\otimes W^*)$, cf. Remark \ref{rem_hamiltonian}. A simple calculation establishes that the corresponding even symmetric bilinear form $\innprod_{\tilde{d}}:=i_2(\tilde{\sigma})$ is given by the formula:
\begin{equation} \label{eqn_dform}
\langle a_1\otimes w_1,a_2\otimes w_2 \rangle_{\tilde{d}} = (-1)^{a_1+w_1}\langle a_1\otimes w_1,\tilde{d}(a_2\otimes w_2) \rangle_{A\otimes W} = (-1)^{w_1 a_2}\langle a_1,a_2 \rangle_d \langle w_1,w_2 \rangle_W.
\end{equation}
Note that this bilinear form is \emph{degenerate}. We have the following lemma.

\begin{lemma}
The quadratic Hamiltonian $\tilde{\sigma}\in S^2(A^*\otimes W^*)$ satisfies both the classical and quantum master equations:
\[ \Delta(\tilde{\sigma})=\{\tilde{\sigma},\tilde{\sigma}\} = 0. \]
\end{lemma}

\begin{proof}
Since $\Delta$ is a \emph{second order} differential operator of \emph{odd} degree and $\tilde{\sigma}$ is an \emph{even quadratic} Hamiltonian, we must have $\Delta(\tilde{\sigma})=0$. It follows from Proposition \ref{prop_hamfields} (cf. also Remark \ref{rem_oddhamfields}) that the condition $\{\tilde{\sigma},\tilde{\sigma}\} = 0$ is equivalent to the condition $\tilde{d}^2=0$.
\end{proof}

As a consequence we have that
\begin{equation} \label{eqn_invmes}
\Delta(e^{-\tilde{\sigma}}) = \left(\frac{1}{2}\{\tilde{\sigma},\tilde{\sigma}\}-\Delta(\tilde{\sigma})\right)\cdot e^{-\tilde{\sigma}} =0.
\end{equation}

Consider the subspace $d(A)\subset A$. Condition \eqref{eqn_dinv} guarantees that it is an isotropic subspace of $A$. Furthermore, the condition that $d$ is contractible implies that the dimension of $d(A)$ is half that of the dimension of $A$; i.e. that $d(A)$ is a \emph{maximally isotropic subspace} of $A$. It follows from Lemma \ref{lem_lagsubspace} \eqref{item_lagsubspace2} that there exists a maximally isotropic subspace $L\subset A$ such that
\begin{equation} \label{eqn_complagsub}
A = L \oplus d(A).
\end{equation}
Moreover, since $d$ must map $L$ isomorphically to $d(A)$, it follows from Lemma \ref{lem_lagsubspace} \eqref{item_lagsubspace3} that the restriction of the bilinear form $\innprod_d$ to $L$ is \emph{nondegenerate}.

From the above data we construct the following functional:

\begin{defi}
Let $A$ be a contractible differential graded commutative Frobenius algebra and let $L\subset A$ be a maximally isotropic subspace of $A$ satisfying \eqref{eqn_complagsub}. We define a functional
\[ S_A^L:C_\bullet(\glie)\to\gf \]
by the formula
\begin{displaymath}
\begin{split}
S_A^L(h_1\wedge\cdots\wedge h_l) &:= (-1)^{p(h)} \langle \Psi_A(h_1)\cdots\Psi_A(h_l)\rangle_0, \\
&:= (-1)^{p(h)} \frac{\int_{L\otimes V} \Psi_A(h_1)\cdots\Psi_A(h_l)\cdot e^{-\tilde{\sigma}}}{\int_{L\otimes V} e^{-\tilde{\sigma}}};
\end{split}
\end{displaymath}
where $h_1,\ldots,h_l\in\glie$ and $V$ is the canonical symplectic vector space of dimension $2n|m$. The sign $(-1)^{p(h)}$ is determined by the Koszul sign rule and appears due to the fact that the map $\Psi_A$ is odd.
\end{defi}

\begin{prop}
The functional $S_A^L:C_\bullet(\glie)\to\gf$ is $\osp$-invariant; that is to say that for any linear symplectic vector field $\eta\in\osp$,
\[ S_A^L(\eta\cdot [h_1\wedge\cdots\wedge h_l]) = 0.\]
Consequently it gives rise to a functional $S_A^L:\cek\to\gf$.
\end{prop}

\begin{proof}
Since $\Psi_A$ is a map of Lie algebras it follows that
\[ S_A^L(\eta\cdot [h_1\wedge\cdots\wedge h_l]) = (-1)^{(l+1)\eta+p(h)} \frac{\int_{L\otimes V} \Psi_A(\eta)[\Psi_A(h_1)\cdots\Psi_A(h_l)]\cdot e^{-\tilde{\sigma}}}{\int_{L\otimes V} e^{-\tilde{\sigma}}} \]
Now by Corollary \ref{cor_leibnizvanishing} \eqref{item_leibnizvanishing2} it follows that
\[ [\Psi_A(\eta)](e^{-\tilde{\sigma}}) = -[\Psi_A(\eta)](\tilde{\sigma})\cdot e^{-\tilde{\sigma}} = 0.\]
Furthermore, by Theorem \ref{thm_vanishing} we have that $\nabla[\Psi_A(\eta)] = 0$. Combining these facts with Lemma \ref{lem_berezin} we get
\begin{displaymath}
\begin{split}
S_A^L(\eta\cdot [h_1\wedge\cdots\wedge h_l]) &= \pm \frac{\int_{L\otimes V} \Psi_A(\eta)[\Psi_A(h_1)\cdots\Psi_A(h_l)\cdot e^{-\tilde{\sigma}}]}{\int_{L\otimes V} e^{-\tilde{\sigma}}}, \\
&= \pm \frac{\int_{L\otimes V} \nabla[\Psi_A(\eta)]\cdot\Psi_A(h_1)\cdots\Psi_A(h_l)\cdot e^{-\tilde{\sigma}}}{\int_{L\otimes V} e^{-\tilde{\sigma}}} = 0.
\end{split}
\end{displaymath}
\end{proof}

We shall now prove that this functional defines a Lie algebra cohomology class.

\begin{theorem} \label{thm_cocycle}
Let $A$ be a contractible differential graded commutative Frobenius algebra:
\begin{enumerate}
\item \label{item_cocycle1}
For any maximally isotropic subspace $L\subset A$ satisfying \eqref{eqn_complagsub}, the functional
\[ S_A^L:\cek\to\gf \]
is a Chevalley-Eilenberg cocycle.
\item \label{item_cocycle2}
The cohomology class of the functional $S_A^L$ does not depend on the choice of maximally isotropic subspace $L$.
\end{enumerate}
\end{theorem}

\begin{proof}
To prove part \eqref{item_cocycle1} we need to show that $S_A^L$ vanishes on Chevalley-Eilenberg boundaries. We do this by first identifying the differential in the Chevalley-Eilenberg complex with the Laplacian defined in section \ref{sec_BVlap} and then applying the corresponding version of Stokes' theorem in the BV-formalism.

Let $\Psi:\cek\to\salg{A^*\otimes V^*}$ be the map given by the formula
\[ \Psi(h_1\wedge\cdots\wedge h_l):= (-1)^{p(h)}\Psi_A(h_1)\cdots \Psi_A(h_l), \]
where $h_1,\ldots, h_l \in \gl{V}$ and the sign $(-1)^{p(h)}$ is determined by the Koszul sign rule as before; then it follows from Theorem \ref{thm_vanishing}, Lemmas \ref{lem_lieids} \eqref{item_lieids2} and \ref{lem_BVidentity} \eqref{item_BVidentity1} and the fact that $\Psi_A$ is a map of Lie algebras that
\begin{equation} \label{eqn_mapcmplx}
\Psi\delta(h_1\wedge\cdots\wedge h_l) = \Delta\Psi(h_1\wedge\cdots\wedge h_l).
\end{equation}
Furthermore, Lemmas \ref{lem_lieids} \eqref{item_lieids2} and \ref{lem_BVidentity} \eqref{item_BVidentity1} also imply that
\begin{equation} \label{eqn_BVCEdiff}
\begin{split}
\Delta[\Psi_A(h_1)\cdots \Psi_A(h_l)\cdot e^{-\tilde{\sigma}}] =& \Delta[\Psi_A(h_1)\cdots \Psi_A(h_l)]\cdot e^{-\tilde{\sigma}} \\
& \pm \Psi_A(h_1)\cdots \Psi_A(h_l)\cdot \Delta(e^{-\tilde{\sigma}}) \\
& + \sum_{i=1}^l \pm \Psi_A(h_1)\cdots\{\Psi_A(h_i),\tilde{\sigma}\}\cdots \Psi_A(h_l)\cdot e^{-\tilde{\sigma}}, \\
=& \Delta[\Psi_A(h_1)\cdots \Psi_A(h_l)]\cdot e^{-\tilde{\sigma}}.
\end{split}
\end{equation}
Here the term on the second line vanishes due to equation \eqref{eqn_invmes} and all the terms in the sum on the third line vanish due to Corollary \ref{cor_leibnizvanishing} \eqref{item_leibnizvanishing2}. It now follows from Corollary \ref{cor_BVstokes} that $S_A^L$ vanishes on Chevalley-Eilenberg boundaries.

We now prove part \eqref{item_cocycle2}. Let $L_0,L_1\subset A$ be two maximally isotropic subspaces satisfying \eqref{eqn_complagsub}. It follows from Lemma \ref{lem_lagsubspace}, parts \eqref{item_lagsubspace1} and \eqref{item_lagsubspace4} that we may choose coordinates $x_1,\ldots,x_n;\xi_1,\ldots,\xi_n$ on $A\otimes V$ such that the Lagrangian subspace $\tilde{L}_0:=L_0\otimes V$ is the locus of the equations
\[ \xi_1=\ldots=\xi_k=0, \quad x_{k+1}=\ldots=x_n=0\]
and such that $\tilde{L}_1:=L_1\otimes V$ is determined by an odd quadratic generating function $\varphi\in S^2(A^*\otimes V^*)$:
\[\xi_1=-\pder{x_1}(\varphi),\ldots,\xi_k=-\pder{x_k}(\varphi),\quad x_{k+1}=\pder{\xi_{k+1}}(\varphi),\ldots, x_n=\pder{\xi_n}(\varphi).\]
Now consider the family of Lagrangian subspaces $\tilde{L}_t$ defined by the generating functions $t\varphi$ for $0\leq t\leq 1$. Using the chain rule, integration by parts, formula \eqref{eqn_canlap} for the Laplacian $\Delta$ and equation \eqref{eqn_invmes}, we obtain the following identity:
\begin{equation} \label{eqn_variation}
\frac{d}{dt}\left[\frac{\int_{\tilde{L}_t}\Psi_A(h_1)\cdots\Psi_A(h_l)\cdot e^{-\tilde{\sigma}}}{\int_{\tilde{L}_t}e^{-\tilde{\sigma}}}\right] = \frac{\int_{\tilde{L}_t}\varphi\cdot\Delta[\Psi_A(h_1)\cdots\Psi_A(h_l)\cdot e^{-\tilde{\sigma}}]}{\int_{\tilde{L}_t}e^{-\tilde{\sigma}}}.
\end{equation}
Consider the functional $G_l:\cek\to\gf$ given by the formula
\[G_l(h_1\wedge\cdots\wedge h_l):=(-1)^{p(h)} \int_0^1 dt \left[\frac{\int_{\tilde{L}_t}\varphi\cdot\Psi_A(h_1)\cdots\Psi_A(h_l)\cdot e^{-\tilde{\sigma}}}{\int_{\tilde{L}_t}e^{-\tilde{\sigma}}}\right].\]
Using equations \eqref{eqn_mapcmplx}, \eqref{eqn_BVCEdiff} and \eqref{eqn_variation} we obtain
\begin{displaymath}
\begin{split}
[S_A^{L_1}-S_A^{L_0}](h_1\wedge\cdots\wedge h_l) &= (-1)^{p(h)} \int_0^1 dt \left[\frac{\int_{\tilde{L}_t}\varphi\cdot\Delta[\Psi_A(h_1)\cdots\Psi_A(h_l)\cdot e^{-\tilde{\sigma}}]}{\int_{\tilde{L}_t}e^{-\tilde{\sigma}}}\right], \\
&=[\delta^*(G_{l-1})](h_1\wedge\cdots\wedge h_l).
\end{split}
\end{displaymath}
\end{proof}

\begin{cor}
To any contractible differential graded commutative Frobenius algebra $A$, there corresponds a well-defined cohomology class $S_A\in H^{\bullet}(\glie,\osp)$.
\end{cor}
\noproof

\subsection{The approach through Feynman diagrams}

In this subsection we recall the original combinatorial version of the dual construction due to Kontsevich \cite{kontfd} and show its equivalence to the version via the BV-formalism. Recall from Definition \ref{def_sympinv} that for any vector space $W$ equipped with an \emph{even} bilinear form, we defined a linear map
\[ \beta_c^W:W^{\otimes 2k}\to\gf \]
for any chord diagram $c\in\chd{k}$ by using the bilinear form to contract the corresponding tensors according to the data contained in the chord diagram.

As in subsection \ref{sec_dualBVform}, we take as our initial data any contractible differential graded commutative Frobenius algebra with an odd symmetric bilinear form, together with a maximally isotropic subspace $L\subset A$ satisfying \eqref{eqn_complagsub}. Consider the even skew-symmetric bilinear form $\innprod_d$ on $A$ given by equation \eqref{eqn_degform}. Its restriction to $L$ is nondegenerate and hence gives rise to an inverse form $\innprod_d^{-1}$ on $L^*$. Let $\tilde{m}_n\in (A^*)^{\otimes n+1}$ be the odd symmetric tensor given by formula \eqref{eqn_itmult}. We can use this data to define a graph cochain by using the tensors $\tilde{m}_n$ as \emph{interactions} and the bilinear form $\innprod_d^{-1}$ as a \emph{propagator}. Given any graph, we place the interaction terms at the vertices and apply the propagator to the edges to produce a number -- the corresponding Feynman amplitude. A precise description is given as follows.

\begin{defi}
Associated to the above data consisting of a contractible differential graded commutative Frobenius algebra $A$  together with a maximally isotropic subspace $L\subset A$, there corresponds a graph cochain
\[ F_A^L:\gc\to\gf\]
defined by the following formula: let $\Gamma$ be any graph, we may assume that there exists a chord diagram $c\in\chd{k}$ and a partition $k_1+\ldots+k_l=2k$ such that $\Gamma=\grpchd{c}$, then
\[ F_A^L(\Gamma):=\beta_c^{L^*}(\tilde{m}_{k_1}\otimes\ldots\otimes\tilde{m}_{k_l}); \]
where the bilinear form $\innprod_d^{-1}$ is used to contract the tensors.
\end{defi}

The following theorem shows that the above construction using Feynman diagrams is equivalent to the previous construction which we defined through the BV-formalism.

\begin{theorem} \label{thm_commute}
The maps in the following diagram commute:
\[\xymatrix{\cek \ar[d]_I \ar[dr]^{S_A^L} \\ \gc \ar[r]_{F_A^L} & \gf}\]
\end{theorem}

\begin{proof}
Here we appeal directly to Wick's formula \eqref{eqn_wick}. Let $h_1,\ldots,h_l\in\glie$ and assume that each $h_i$ is a monomial of degree $k_i\geq 3$.
\begin{displaymath}
\begin{split}
S_A^L(h_1\wedge\cdots\wedge h_l) &= (-1)^{p(h)}\langle \Psi_A(h_1)\cdots\Psi_A(h_l) \rangle_0 \\
&= (-1)^{p(h)}\sum_{c\in\chd{k}}\beta_c^{L^*\otimes V^*}[\Psi_A(h_1)\otimes\ldots\otimes\Psi_A(h_l)] \\
&= \sum_{c\in\chd{k}}\beta_c^{L^*}[\tilde{m}_{k_1}\otimes\ldots\otimes \tilde{m}_{k_l}]\beta_c^{V^*}[h_1\otimes\ldots\otimes h_l] \\
&= F_A^L I(h_1\wedge\cdots\wedge h_l).
\end{split}
\end{displaymath}
The second line follows from Lemma \ref{lem_wicks}. The third line follows from equation \eqref{eqn_dform}.
\end{proof}

\appendix
\section{Superintegral calculus}

In this section we briefly review the basic definition of superintegral and prove an auxiliary lemma. Here we will only be concerned with integrating polynomial functions with respect to Gaussian measures, although the theory could be described in more general contexts.

Let $\mathcal{F}_{n|m}$ denote the space of functions on $\gf^{n|m}$ of the form $f(x,\xi)e^{-\frac{1}{2}\sum_{i=1}^n x_i^2}$ where $f$ is a polynomial superfunction on $\gf^{n|m}$. The integral operators
\[ \int_{-\infty}^{\infty}dx_i:\mathcal{F}_{n|m}\to\mathcal{F}_{n-1|m} \]
of integration with respect to an \emph{even} variable have their standard meaning, which extends to superfunctions in an obvious manner. The integral operators
\[ \int_{-\infty}^{\infty}d\xi_i:\mathcal{F}_{n|m}\to\mathcal{F}_{n|m-1} \]
of integration with respect to an \emph{odd} variable are defined simply as $\int_{-\infty}^{\infty}d\xi_i:=\pder{\xi_i}$, where $\pder{\xi_i}$ by convention acts on the \emph{right}. This gives us the following definition:

\begin{defi}
Let $g\in\mathcal{F}_{n|m}$, the superintegral of $g$ is defined by the formula
\[ \int_{\mathbb{R}^{n|m}} g\, dx_1\ldots dx_n d\xi_1\ldots d\xi_m := \int_{-\infty}^{\infty}\cdots\int_{-\infty}^{\infty} g\, dx_1\ldots dx_n d\xi_1\ldots d\xi_m.\]
\end{defi}

We now prove the infinitesimal analogue of the Berezin change of variables formula.

\begin{lemma} \label{lem_berezin}
Let $g\in\mathcal{F}_{n|m}$ be any function and $\eta$ be a polynomial vector field on $\gf^{n|m}$, then
\[ \int_{\mathbb{R}^{n|m}} \eta(g)\, dxd\xi = - \int_{\mathbb{R}^{n|m}} \nabla(\eta)g\, dxd\xi.\]
\end{lemma}

\begin{proof}
We may assume that $\eta$ has either the form $\eta:=p\pder{\xi_i}$ or the form $\eta:=p\pder{x_i}$ for some polynomial $p$. We claim that
\[ \int_{\gf^{n|m}} \nabla(g\cdot\eta)\,dxd\xi = 0. \]
In the former case this is a simple consequence of the definitions. In the latter case it follows from Lemma \ref{lem_stokes}. The result now follows as a simple consequence of Lemma \ref{lem_dividentities} \eqref{item_dividentitiesfmult}.
\end{proof}

\end{document}